\documentclass[reqno,12pt]{amsart} 
\numberwithin{equation}{section}
\usepackage{amssymb}
\usepackage[mathscr]{eucal}
\usepackage{amsmath}
\usepackage{latexsym}
\usepackage{amsthm}
\theoremstyle{plain} 
\newtheorem{theorem}{\indent\sc Theorem}[section]
\newtheorem{lemma}[theorem]{\indent\sc Lemma}

\newtheorem{proposition}[theorem]{\indent\sc Proposition}

\theoremstyle{definition} 

\begin{document}

\title[Weighted norm inequalities for fractional oscillatory integrals]
{Weighted norm inequalities for fractional oscillatory integrals and applications}

\author[Shaoguang Shi]{Shaoguang Shi} 

\author[Zunwei Fu]{Zunwei Fu}

\author[Shanzhen Lu]{Shanzhen Lu}

\author[Fayou Zhao]{Fayou Zhao}
%
\thanks{ 
$^{*}$This work was partially supported by
 NSF of China (Grant No. 10901076,  10931001 and 11171345), the Key Laboratory of Mathematics and Complex System (Beijing Normal
University), Ministry of Education, China.}
\address{
School of Mathematical Sciences\endgraf
Beijing Normal
University\endgraf
Beijing 100875\endgraf
and\endgraf
\endgraf
School of Sciences\endgraf
Linyi
University \endgraf
Linyi 276005\endgraf
P. R. China}
\email{shishaoguang@lyu.edu.cn}

\address{
School of Sciences\endgraf
Linyi
University \endgraf
Linyi 276005\endgraf
P. R. China}
\email{lyfzw@tom.com}

\address{
School of Mathematical Sciences\endgraf
Beijing Normal
University\endgraf
Beijing 100875\endgraf
\endgraf
}
\email{lusz@bnu.edu.cn}

\address{
College of Sciences\endgraf
Shanghai
University\endgraf
Shanghai 200444\endgraf
\endgraf
}
\email{zhaofayou2008@yahoo.com.cn}

\maketitle

\begin{abstract}
We set up some weighted norm inequalities for fractional oscillatory
integral operators.
As applications, the corresponding results for commutators formed
by $BMO(\mathbb{R}^{n})$ functions and the operators are established.
\end{abstract}

\section{Introduction} 
We consider a class of fractional oscillatory integrals
$$(T_{\alpha}f)(x)=\int_{\mathbb{R}^{n}}e^{iP(x,y)}K_{\alpha}(x,y)f(y)dy,\eqno(1.1)$$
where $P(x,y)$ is a real valued polynomial defined on $\mathbb{R}^{n}\times\mathbb{R}^{n}$, and $K_{\alpha}$ is a fractional Calder\'{o}n-Zygmund kernel which satisfies
$$|K_{\alpha}(x,y)|\leq\frac{C}{|x-y|^{n-\alpha}}, \,\,\,\,\,\,0<\alpha<n \eqno(1.2)$$
and
$$|\nabla_{x}K_{\alpha}(x,y)|+|\nabla_{y}K_{\alpha}(x,y)|\leq\frac{C}{|x-y|^{n+1-\alpha}}.\eqno(1.3)$$
In the foregoing and following, the letter $C$ will stand for a positive constant which may vary from line to line.

In 1987, Ricci and Stein \cite{RS} obtained the $L^{p}(1<p<\infty)$-boundedness of $T_{\alpha}$, with bound depending on the polynomial. For the $L^{p}(1<p<\infty)$-boundedness of $T_{\alpha}$ with rough kernels, we refer the reader to Ding's work in \cite{D1} and \cite{D2}.
Obviously, when $\alpha=0$, $K_{0}=K$ is exactly the classical standard Calder\'{o}n-Zygmund kernel and the corresponding integral
$$(Tf)(x)=\mathrm{p.v.}\int_{\mathbb{R}^{n}}e^{iP(x,y)}K(x,y)f(y)dy,\eqno(1.4)$$
is the classical oscillatory singular integral. Moreover, it is well-known that the Radon transform \cite{St}, being an important role in the
CT technology of medical sciences, is closely related to this form of oscillatory
singular integrals. We may refer to Ricci and Stein's work \cite{RS} for their $L^{p}(1<p<\infty)$-norm estimates and Lu-Zhang's work \cite{LZ2} about the weighted versions, respectively. For the other works about (1.4), we would like to refer to \cite{CC}, \cite{L1}, \cite{L2}, \cite{LZ1}, \cite{Sa} and \cite{St}. Highly inspired by \cite{LZ2} and \cite{RS}, our first goal is to get the weighted strong type-boundedness for $T_{\alpha}$. Before stating our results, let us give some notations first.

Let $1<p<\infty$. For any non-negative locally functions $w$ and any Lebesgue measurable function $f$, we set
$$\|f\|_{L^{p}(w)}=\left(\int_{\mathbb{R}^{n}}|f(x)|^{p}w(x) dx\right)^{1/p}$$
and if $w\equiv1$, we denote $\|f\|_{L^{p}(w)}$ simply by $\|f\|_{L^{p}}$. The Muckenhoupt classes $A_{p}$ and $A_{(p,q)}$ \cite{St} contain the functions $w$ which satisfy
$$A_{p}:\sup_{Q}\left(\frac{1}{|Q|}\int_{Q}w(x)dx\right)\left(\frac{1}{|Q|}\int_{Q}w(x)^{1-p'}dx\right)^{p-1}\leq C, 1<p<\infty$$
and
$$
A_{(p,q)}:\sup_{Q}\left(\frac{1}{|Q|}\int_{Q}w(x)^{q}dx\right)^{\frac{1}{q}}\left(\frac{1}{|Q|}\int_{Q}w(x)^{-p'}dx\right)^{\frac{1}{p'}}\leq C, 1<p, q<\infty,
$$
respectively. Here $Q$ denotes any cubes in $\mathbb{R}^{n}$, $1/p+1/p'=1$. Define the Hardy-Littlewood maximal operator $M$ as
$$Mf(x)=\sup_{Q\ni x}\frac{1}{|Q|}\int_{Q}|f(y)|dy.$$
It is well known \cite{St} that $M$ is a bounded operator on $L^{p}({\mathbb
R}^{n})$, $1<p<\infty$.

The assumption that $P(x,y)$ is nontrivial will be needed throughout the paper. We say a polynomial $P(x,y)$ is nontrivial if $P(x,y)$ does not take the form $P_{0}(x)+P_{1}(y)$, where $P_{0}$ and $P_{1}$ are polynomials defined on $\mathbb{R}^{n}$. Now we may formulate our results as follows
\begin{theorem} Suppose $P(x,y)=\sum_{|\beta|\leq k,|\gamma|\leq l}a_{\beta\gamma}x^{\beta}y^{\gamma}$ is nontrivial real valued polynomial, $K_{\alpha}$ satisfies $(1.2)$, $(1.3)$ and the operator
$\widetilde{T}_{\alpha}: f\rightarrow \int_{\mathbb{R}^{n}}K_{\alpha}(x,y)f(y)dy$ is bounded on $L^{2}$. Then, if $w\in A_{p}$, $1<p<\infty$, there exists constant $C(\alpha)>0, 0<\delta<1$, such that
$$\|T_{\alpha}f\|_{L^{p}(w)}\leq C(\alpha)\|f\|_{L^{p}(w)},$$
where $0<\alpha<(1+\frac{l}{k})\delta$.
\end{theorem}

It would be desirable to see what happens if we replace $w\in A_{p}$ by $w\in A_{(p,q)}$, that is
\begin{theorem} Let $K_{\alpha}(x,y),P(x,y)$ be the same as Theorem $1.1$, $w\in A_{(p,q)}$ and $0<\alpha<n$, $1< p<\frac{n}{\alpha}, \frac{1}{p}-\frac{1}{q}=\frac{\alpha}{n}$. Then
there exists constant $C>0$ independent of $\alpha$
such that
$$\|T_{\alpha}f\|_{L^{q}(w^{q})}\leq C\|f\|_{L^{p}(w^{p})}.$$
\end{theorem}
It is well known that the oscillatory factor $e^{ip(x,y)}$ makes it impossible to establish the weighted norm inequalities of (1.1) by the methods as in the case of  Calder\'{o}n-Zygmund  operators or fractional integrals. It is also worth pointing out that the methods used in \cite{D1} and \cite{D2} depends heavily on the homogeneity of the rough kernel. Since there is not any homogeneity for $K_{\alpha}$, the method of \cite{D1} and \cite{D2} does not work in dealing with the kernel in the proof of Theorem 1.1. Using properties of weights functions and interpolation idea due to Stein and Weiss \cite{SW}, as well as properties established by us, we first give the proof of Theorem 1.1 and Theorem 1.2 in Section 2. As applications of Theorem 1.1 and Theorem 1.2, we shall discuss the corresponding results for commutators formed by $BMO(\mathbb{R}^{n})$ functions and fractional oscillatory integrals in Section 3.

\section{Proofs of Theorem 1.1 and Theorem 1.2}
Our treatment of the fractional oscillatory operator $T_{\alpha}$ will be based in part on some simple but useful inequalities concerning polynomials. For the convenience of the reader, we repeat the relevant material from \cite{RS} and \cite{St} without proofs, thus making our exposition self-contained.
Let $P(x)=\sum_{|\beta|\leq d}a_{\beta}x^{\beta}$ denote a polynomial in $\mathbb{R}^{n}$ of degree $d$, where $x^{\beta}=x_{1}^{\beta_{1}}x_{2}^{\beta_{2}}\cdots x_{n}^{\beta_{n}},$ $\beta=(\beta_{1},\cdots, \beta_{n})$ with $|\beta|=\beta_{1}+ \beta_{2}+\cdots+\beta_{n}$.

\begin{lemma}
Suppose that $\varepsilon< \frac{1}{d}$. Then
$$
\int_{|x|\leq 1}|P(x)|^{-\varepsilon}dx\leq A_{\varepsilon}\left(\sum_{|\beta|\leq d}|a_{\beta}|\right)^{-\varepsilon}.
$$
The bound $A_{\varepsilon}$ depends on $\varepsilon$ (and the dimension $n$) but not on the coefficients $\{a_{\beta}\}$.
\end{lemma}
The above lemma will be applied via following consequence.
\begin{lemma}
$\mathrm{(1)}$ \,\,Let $P(x)=\sum_{|\beta|= d}a_{\beta}x^{\beta}$ denote a homogeneous polynomial in $\mathbb{R}^{n}$ of degree $d$ with $\varepsilon< \frac{1}{d}$. Then
$$
\int_{|x|= 1}|P(x)|^{-\varepsilon}d\sigma (x)\leq A_{\varepsilon}\left(\sum_{|\beta|= d}|a_{\beta}|\right)^{-\varepsilon}.\eqno(2.1)
$$
$\mathrm{(2)}$ \,\,Let $P(x)=\sum_{|\beta|\leq d}a_{\beta}x^{\beta}$ denote a polynomial in $\mathbb{R}^{n}$ of degree $d$. Then
$$
\sup_{y\in \mathbb{R}^{n}}\int_{|x|\leq 1}|P(x-y)|^{-\varepsilon}dx\leq A_{\varepsilon}\left(\sum_{|\beta|= d}|a_{\beta}|\right)^{-\varepsilon}.\eqno(2.2)
$$

\end{lemma}
Besides the above inequalities for polynomials. We still need the following estimates for one-dimensional oscillatory integral whose phases are generalized polynomials.

\begin{lemma}
Let $\psi\in C^{1}[\alpha,\beta]$, $\varepsilon=\min\{\frac{1}{a_{1}},\frac{1}{n}\}$, $\lambda>0$. Then
$$
\left|\int_{\alpha}^{\beta}e^{i\lambda\phi(t)}\psi(t)dt\right|\leq C\lambda^{-\varepsilon}\left\{\sup_{\alpha\leq t\leq\beta}|\psi(t)|+\int_{\alpha}^{\beta}|\psi'(t)|dt\right\},
$$
where $\phi$ is real-valued phases of the form $\phi(t)=t^{a_{1}}+\mu_{2}t^{a_{2}}+\cdots+\mu_{n}t^{a_{n}}$ with $\mu_{2},\cdots,\mu_{n}$ are real parameters and $a_{1},a_{2},\cdots,a_{n}$ are distinct positive exponents.
\end{lemma}
Let $0<\alpha<n$. The fractional integrals were defined by
$$I_{\alpha}f(x)=\int_{\mathbb{R}^{n}}\frac{f(y)}{|y-x|^{n-\alpha}}dy.$$
Muchenhout and Wheeden established the weighted boundedness of $I_{\alpha}$ in \cite{MW}.
\begin{lemma}
Let $0<\alpha<n$, $1< p<\frac{n}{\alpha}<\infty, \frac{1}{p}-\frac{1}{q}=\frac{\alpha}{n}$ and $w\in A_{(p,q)}$. Then there exist constant $C>0$, such that if $1<p<\infty$, then
$$\|I_{\alpha}f\|_{L^{q}(w^{q})}\leq C\|f\|_{L^{p}(w^{p})}.$$
\end{lemma}

Recalling the definition of $A_{p}$ and $A_{(p,q)}$, we proceed to show some relationship between the classes $A_{p}$ and $A_{(p,q)}$.
\begin{lemma}\,\,\,\,\cite{LDY}\,
Suppose $0<\alpha<n,$ $1<p<\frac{n}{\alpha}<\infty$, and $\frac{1}{p}-\frac{1}{q}=\frac{\alpha}{n}$. Then

$\mathrm{(1)}$\,\, If $w\in A_{(p,q)}\Longleftrightarrow w^{q}\in A_{q(\frac{n-\alpha}{n})}\Longleftrightarrow w^{-p'}\in A_{1+\frac{p'}{q}}$.

$\mathrm{(2)}$\,\, $w^{q}\in A_{(p,q)}\Longrightarrow w^{q}\in A_{q}$ and $w^{p}\in A_{p}$.

\end{lemma}
From the reverse H\"{o}lder inequality and Lemma 2.5, we can easily have the following properties of $A_{p}$ and $A_{(p,q)}$ weights, which will paly an important role in the proof of Theorem 1.1 and Theorem 1.2.
\begin{lemma}\,\,\,\,\cite{LDY}\,
$\mathrm{(1)}$\,\, If $w\in A_{p}$, there exists some $\varepsilon>0$ such that $w^{1+\varepsilon}\in A_{p}$.

$\mathrm{(2)}$\,\, Suppose $0<\alpha<n,$ $1<p<\frac{n}{\alpha}<\infty$, and $\frac{1}{p}-\frac{1}{q}=\frac{\alpha}{n}$. Then if $w\in A_{(p,q)}$, there exists some $\varepsilon>0$ such that $w^{1+\varepsilon}\in A_{(p,q)}$.
\end{lemma}
To prove the main results, the interpolation theorem of operators with change measures plays an important role, which formed by Stein and Weiss \cite{SW}.
\begin{lemma}
Suppose that $u_{0},v_{0},u_{1},v_{1}$ are positive weight functions and $1<p_{0},p_{1}<\infty$. Assume sublinear operator $S$ satisfies:
$$\|Sf\|_{L^{p_{0}}(u_{0})}\leq C_{0}\|f\|_{L^{p_{0}}(v_{0})},$$
and
$$\|Sf\|_{L^{p_{1}}(u_{1})}\leq C_{1}\|f\|_{L^{p_{1}}(v_{1})}.$$
Then
$$\|Sf\|_{L^{p}(u)}\leq C\|f\|_{L^{p}(v)}$$
holds for any $0<\theta<1$ and $\frac{1}{p}=\frac{\theta}{p_{0}}+\frac{1-\theta}{p_{1}}$, where
$u=u_{0}^{\frac{p\theta}{p_{0}}}u_{1}^{\frac{p(1-\theta)}{p_{1}}}$, $v=v_{0}^{\frac{p\theta}{p_{0}}}v_{1}^{\frac{p(1-\theta)}{p_{1}}}$ and $C\leq C_{0}^{\theta}C_{1}^{1-\theta}$.
\end{lemma}

We now give the proof of Theorem 1.1 and Theorem 1.2 in what follows.

{\it Proof of Theorem $1.1$}. Suppose $P(x,y)$ is a nontrivial real polynomial with degree $k$ in $x$ and degree $l$ in $y$. We shall carry out the argument by induction. First, we assume the conclusion of Theorem 1.1 is valid for all polynomials which are the sums of monomials of degree less than $k$ in $x$ times
monomials of any degree in $y$, together with monomials which are of degree $k$ in $x$  times monomials which are of degree less than $l$ in $y$ .Thus $P(x,y)$ can be written as
$$
P(x,y)=\sum_{|\beta|=k,|\gamma|=l}a_{\beta\gamma}x^{\beta}y^{\gamma}+R(x,y)
$$
or
$$
P(x,y)=\sum_{|\beta|=k}x^{\beta}Q_{\beta}(y)+R_{0}(x,y).
$$
Here $R(x,y)$ satisfies the above induction assumption and $Q_{\beta}(y)=\sum_{|\gamma|<l}a_{\beta\gamma}y^{\gamma}$ is polynomial in $y$ of degree less than $l$. $R_{0}(x,y)$ has $x$-degree less than $k$. Without loss of generality, we may assume $k,l>0$ and $\sum_{|\beta|=k,|\gamma|=l}|a_{\beta\gamma}|=1.$ We split the kernel $K_{\alpha}$ as
$$
K_{\alpha}(x,y)=K_{\alpha}(x,y)\chi_{\{|x-y|\leq1\}}(x,y)+K_{\alpha}(x,y)\chi_{\{|x-y|>1\}}(x,y)=:K_{\alpha,0}+K_{\alpha,\infty},
$$
and consider the corresponding splitting
$$
T_{\alpha}f(x)=\int e^{iP(x,y)}K_{\alpha,0}(x,y)f(y)dy+\int e^{iP(x,y)}K_{\alpha,\infty}(x,y)f(y)\,dy
$$
$$
=:T_{\alpha,0}f(x)+T_{\alpha,\infty}f(x).
$$
Our task is now to show the weighted estimates for both $T_{\alpha,0}$ and $T_{\alpha,\infty}$.
Take any $h\in \mathbb{R}$, and write
$$
P(x,y)=\sum_{|\beta|=k,|\gamma|=l}a_{\beta\gamma}(x-h)^{\beta}(y-h)^{\gamma}+R(x,y,h),
$$
where the polynomial $R(x,y,h)$ satisfies the induction assumption, and the coefficients of $R(x,y,h)$ depend on $h$.

We first set up the estimate for $T_{\alpha,0}$.

Observing $$e^{iP(x,y)}=e^{iP(x,y)}-e^{i(R(x,y,h)+\Sigma a_{\beta\gamma}(y-h)^{\beta+\gamma})}+e^{i(R(x,y,h)+\Sigma a_{\beta\gamma}(y-h)^{\beta+\gamma})},$$ we have
$$\begin{array}{rl}
\displaystyle T_{\alpha,0}f(x)&=\displaystyle \int_{|x-y|<1}e^{i(R(x,y,h)+\Sigma a_{\beta\gamma}(y-h)^{\beta+\gamma})}K_{\alpha}(x,y)f(y)dy\\
&\,\,+\displaystyle \int_{|x-y|<1}\left\{e^{iP(x,y)}-e^{i(R(x,y,h)+\Sigma a_{\beta\gamma}(y-h)^{\beta+\gamma})}\right\}K_{\alpha}(x,y)f(y)dy\\&=:\displaystyle T_{\alpha,01}f(x)+T_{\alpha,02}f(x).
\end{array}$$
Now we split $f$ into three parts as follows
$$\begin{array}{rl}
\displaystyle f(y)&=\displaystyle f(\cdot)\chi_{\{|y-h|<\frac{1}{2}\}}(y)+f(\cdot)\chi_{\{\frac{1}{2}\leq|y-h|<\frac{5}{4}\}}(y)+f(\cdot)\chi_{\{|y-h|\geq\frac{5}{4}\}}(y)\\
&=:f_{1}(y)+f_{2}(y)+f_{3}(y).
\end{array}$$

It is easy to see that when $|x-h|<\frac{1}{4}$, we have
$$
T_{\alpha,01}f_{1}(x)=\int e^{i(R(x,y,h)+\Sigma a_{\beta\gamma}(y-h)^{\beta+\gamma})}K_{\alpha}(x,y)f_{1}(y)dy.
$$
Thus, it follows from the induction assumption that
$$
\int_{|x-h|<\frac{1}{4}}|T_{\alpha,01}f_{1}(x)|^{p}w(x)dx\leq C\int_{|y-h|<\frac{1}{2}} |f(y)|^{p}w(y)dy,\eqno(2.3)
$$
where $C$ is independent of $h$.

Notice that $|x-h|<\frac{1}{4},\frac{1}{2}\leq|y-h|<\frac{5}{4}$ imply $|y-x|>\frac{1}{4}$. Thus
$$
|T_{\alpha,01}f_{2}(x)|\leq C\int_{\frac{1}{4}<|x-y|<1}|K_{\alpha}(x,y)f_{2}(y)|dy\leq C(\alpha)M(f_{2})(x).
$$
Here $M$ denotes the Hardy-Littlewood maximal operator. So by the weighted boundedness of  $M$, we have
$$
\int_{|x-h|<\frac{1}{4}}|T_{\alpha,01}f_{2}(x)|^{p}w(x)dx\leq C\int_{|y-h|<\frac{5}{4}} |f(y)|^{p}w(y)dy,\eqno(2.4)
$$
where $C$ is independent of $h$.

Finally, noticing that if $|x-h|<\frac{1}{4},|y-h|\geq\frac{5}{4}$, we have $|y-x|>1$,thus
$$
T_{\alpha,01}f_{3}(x)=0.\eqno(2.5)
$$

Combining (2.3), (2.4) and (2.5), we get
$$
\int_{|x-h|<\frac{1}{4}}|T_{\alpha,01}f(x)|^{p}w(x)dx\leq C\int_{|y-h|<\frac{5}{4}} |f(y)|^{p}w(y)dy,\eqno(2.6)
$$
where $C$ is independent of $h$.

If $|x-h|<\frac{1}{4}, |x-y|<1$, then
$$
|e^{iP(x,y)}-e^{i(R(x,y,h)+\sum_{\beta\gamma}a_{\beta\gamma}(y-h)^{\beta+\gamma})}|\leq \sum|a_{\beta\gamma}||x-y|=C|x-y|.
$$
We have
$$
|T_{\alpha,02}f(x)|\leq C\int_{|x-y|<1}\frac{|f(y)|}{|x-y|^{n-\alpha-1}}dx\leq CM(f(\cdot)\chi_{B(h,\frac{5}{4})}(\cdot))(x),
$$
which implies
$$
\int_{|x-h|<\frac{1}{4}}|T_{\alpha,02}f(x)|^{p}w(x)dx\leq C\int_{|y-h|<\frac{5}{4}} |f(y)|^{p}w(y)dy,\eqno(2.7)
$$
where $C$ is independent of $h$.

From (2.6) and (2.7), it follows that the inequality
$$
\int_{|x-h|<\frac{1}{4}}|T_{\alpha,0}f(x)|^{p}w(x)dx\leq C\int_{|y-h|<\frac{5}{4}} |f(y)|^{p}w(y)dy,
$$
holds uniformly in $h\in \mathbb{R}$, which means
$$
\|T_{\alpha,0}f\|_{L^{p}(w)}\leq C\|f\|_{L^{p}(w)},\eqno(2.8)
$$
where $w\in A_{p}$.

We now turn to the estimate for $T_{\alpha, \infty}$.

We write $K_{\infty}(x,y)=\sum_{j=0}^{\infty}\psi_{j}(x,y)$, where $\psi_{0}$ is supported in $\frac{1}{2}\leq |x-y|\leq 1$ and is bounded, while $\psi_{j}(x,y)=2^{-j(n-\alpha)}\psi(\frac{x}{2^{j}},\frac{y}{2^{j}})$, $j\geq 1$ with $\psi$ an appropriate function of class $C^{1}$ supported in $\frac{1}{2}\leq |x|\leq 1$. We set
$$T_{\alpha,j}f(x)=\int e^{iP(x,y)}\psi_{j}(x,y)f(y)dy.$$
Since bounds for $T_{\alpha,j}$ when $j=0$ is trivial, we turn to the case for $j\geq 1.$
In order to apply the method of interpolation of operators with change of measures, we first set up the $L^{2}-$ boundedness for $T_{\alpha,j}$, which is also essential to set up the $L^{p}-$boundedness$(1<p<\infty)$ for $T_{\alpha, \infty}$.
\begin{proposition}
There exists constant $C>0$, such that
$$\|T_{\alpha,j}\|_{L^{2}}\leq C2^{-j(\frac{1-\alpha}{2}+\frac{l}{2k})}.$$
\end{proposition}
{\it Proof.}
Since $|T_{\alpha,j}|^{2}=T_{\alpha,j}T^{\ast}_{\alpha,j}$, we only need to consider $T_{\alpha,j}T^{\ast}_{\alpha,j}$, which kernel is given by
$$K_{j}(y,z)=\int e^{i(P(x,z)-P(x,y))}\psi_{j}(x,z)\overline{\psi_{j}}(x,y)dx.$$
By rescaling, we would only to prove
$$\sup_{z,y}\int |\widetilde{K_{j}}(y,z)|dzdy\leq C2^{-j(1+\frac{l}{k}-\alpha)},$$
where
$$\widetilde{K_{j}}(y,z)=2^{(n-\alpha)j}K_{j}(2^{j}y,2^{j}z)=2^{\alpha j}\int e^{i(P(2^{j}x,2^{j}z)-P(2^{j}x,2^{j}y))}\psi_{0}(x,z)\overline{\psi_{0}}(x,y)dx.\eqno(2.9)$$
A trivial estimates for $\widetilde{K_{j}}$ is $|\widetilde{K_{j}}(y,z)|\leq C\chi_{\{|y-z|<2\}}(y,z)$. We make the changes of variables $x\rightarrow x+y$, which yields
$$P(2^{j}(x+y),2^{j}z)=\sum_{|\beta|=k}2^{|\beta|j}x^{\beta}Q_{\beta}(2^{j}z)+R_{j}(x,y,z),\eqno(2.10)$$
where $R_{j}$ has $x-$degree strictly less than $k$. Similarly
$$P(2^{j}(x+y),2^{j}y)=\sum_{|\beta|=k}2^{|\beta|j}x^{\beta}Q_{\beta}(2^{j}y)+R'_{j}(x,y).\eqno(2.11)$$
Substituting (2.10) and (2.11) into (2.9) yields
$$\widetilde{K_{j}}(y,z)=2^{\alpha j}\int e^{i(\sum_{|\beta|=k}2^{|\beta|j}x^{\beta}(Q_{\beta}(2^{j}z)-Q_{\beta}(2^{j}y))+R_{j}(x,y,z)-R'_{j}(x,y))}\psi_{0}(x+y,z)\overline{\psi_{0}}(x+y,y)dx.$$
Applying the polar coordinates with $x=rx',r=|x|, |x'|=1$, we can rewrite the above equality as
$$\widetilde{K_{j}}(y,z)=2^{\alpha j}\int_{|x'|=1}\left(\int_{0}^{1} e^{i(A+B)}\widetilde{\psi}dr\right)d\delta x',\eqno(2.12)$$
where
$$A=(r2^{j})^{k}\sum_{|\beta|=k}x'^{\beta}(Q_{\beta}(2^{j}z)-Q_{\beta}(2^{j}y)),$$
$$B=R_{j}(x,y,z)-R'_{j}(x,y),$$
and
$$\widetilde{\psi}=\psi_{0}(x+y,z)\overline{\psi_{0}}(x+y,y)r^{n-1}.$$

Notice that $\widetilde{R}$ has degree strictly less than $k$ when viewed as a polynomial in $r$. An application of Lemma 2.3 with $a_{1}=n=k, \varepsilon=\frac{1}{k}$ for the inner integral in (2.12) we obtain
\begin{eqnarray*}
|\widetilde{K_{j}}(y,z)|&\leq& C2^{\alpha j}2^{-j}\int_{|x'|=1}\sum_{|\beta|=k}x'^{\beta}|Q_{\beta}(2^{j}z)-Q_{\beta}(2^{j}y)|^{-\frac{1}{k}}dx'\\
&\leq& C2^{\alpha j}2^{-j}\left(\sum_{|\beta|=k}|Q_{\beta}(2^{j}z)-Q_{\beta}(2^{j}y)|\right)^{-\frac{1}{k}}\chi_{\{|y-z|<2\}}(y,z),
\end{eqnarray*}
where we use the estimate (2.1) in the last inequality. Since $\sum_{|\beta|=k,|\gamma|=l}|a_{\beta\gamma}|=1$, there is an $|\beta_{0}|=k$ and $|\gamma_{0}|=l$ such that $|a_{\beta_{0}\gamma_{0}}|\geq l>0.$ Recalling the definition of $Q_{\beta_{0}}(2^{j}z)=2^{jl}\sum_{|\gamma|=l}a_{\beta_{0}\gamma}z^{\gamma}$+lower order terms in $z$, we have
$$
\int|\widetilde{K_{j}}(y,z)|dz\leq C2^{\alpha j}2^{-j}\left(\int_{|z-y|\leq 2}\int_{|\beta|=k}\left(|Q_{\beta_{0}}(2^{j}z)-Q_{\beta_{0}}(2^{j}y)|\right)^{-\frac{1}{k}}\right)dz.
$$
By (2.2), we obtain
$$\sup_{y}\int|\widetilde{K_{j}}(y,z)|dz\leq C2^{-j(1+\frac{l}{k})}2^{\alpha j}=C2^{-j(1+\frac{l}{k}-\alpha)}.$$
Similarly,
$$\sup_{z}\int|\widetilde{K_{j}}(y,z)|dy\leq C2^{-j(1+\frac{l}{k})}2^{\alpha j}=C2^{-j(1+\frac{l}{k}-\alpha)},$$
which complete the proof of Proposition 2.8.
 \qed

Obviously, the norms of $T_{\alpha,j}$ on $L^{1}$ and $L^{\infty}$ are both $2^{\alpha j}$. Thus by interpolation, we have for $1<p<\infty, 0<\theta<1$
$$\|T_{\alpha,j}\|_{L^{p}}\leq C2^{-j\left[\theta(\frac{1-\alpha}{2}+\frac{l}{2k})-\alpha(1-\theta)\right]}.\eqno(2.13)$$
Now, we proceed to the estimates for $T_{\alpha,\infty}$. For $j\geq 1$, we have
$$
|T_{\alpha,j}f(x)|\leq \int_{2^{j-1}<|x-y|<2^{j}}\frac{|f(y)|}{|x-y|^{n-\alpha}}dy\leq C2^{j\alpha}M(f)(x).
$$
Thus from Lemma 2.6 (1), we have
$$
\|T_{\alpha,j}f\|_{L^{p}(w^{1+\varepsilon})}\leq C2^{j\alpha}\|f\|_{L^{p}(w^{1+\varepsilon})},\eqno(2.14)
$$
where $C$ is independent of $j$.

Applying Lemma 2.7 to (2.13) and (2.14), we have
$$
\|T_{\alpha,j}f\|_{L^{p}(w)}\leq C2^{-j[\frac{(1+\frac{l}{k}-\alpha )}{2})\theta^{2}-\alpha(1-\theta^{2})]}\|f\|_{L^{p}(w)}.
$$
Taking $\delta=\frac{\theta^{2}}{2-\theta^{2}}$, for $0<\alpha<(1+\frac{l}{k})\delta$, we have
$$
\|T_{\alpha,\infty}\|_{L^{p}(w)}\leq C.
$$
 So, we complete the proof of Theorem 1.1. \qed

{\it Proof of Theorem $1.2$}. Similar to that of Theorem 1.1, we can get Theorem 1.2. However, we can prove Theorem 1.2 in a simpler manner. In fact, Theorem 1.2 is a straightforward results by Lemma 2.4 and the following observation
$$
|T_{\alpha}f(x)|\leq  \int\frac{|f(y)|}{|y-x|^{n-\alpha}}dy= I_{\alpha}(|f|)(x).
$$
 \qed

\section{Weighted estimates for the commutators of fractional oscillatory integrals}
As the applications of Theorem 1.1 and Theorem 1.2, we show the weighted boundedness of commutators of fractional oscillatory integrals in this section.

Let $b\in L_{loc}^{1}(\mathbb{R}^{n})$, we say that $b\in BMO(\mathbb{R}^{n})$ if
$$\|b\|_{BMO(\mathbb{R}^{n})}=\sup_{Q}\frac{1}{|Q|}\int_{Q}|b-b_{Q}|<\infty.$$
Here $b_{Q}=\frac{1}{|Q|}\int_{Q}b$. There were some relationship between the Muchenhoupt classes and $BMO(\mathbb{R}^{n})$ functions \cite{LDY}.

\begin{lemma}\,\,\,\,

$\mathrm{(1)}$\,\,If $w\in A_{p},1<p<\infty$, $b\in BMO(\mathbb{R}^{n})$, then for $\lambda>0$, there exists $\eta>0$ such that
$$e^{\lambda b}\in A_{p}$$
when $\|b\|_{BMO(\mathbb{R}^{n})}<\eta$.

$\mathrm{(2)}$\,\,Suppose $0<\alpha<n,$ $1<p<\frac{n}{\alpha}<\infty$, and $\frac{1}{p}-\frac{1}{q}=\frac{\alpha}{n}$. Then for $w\in A_{(p,q)}$,  $\lambda>0$, there exists $\eta>0$ such that
$$e^{\lambda b}\in A_{(p,q)}$$
when $\|b\|_{BMO(\mathbb{R}^{n})}<\eta$.\end{lemma}

Let $b\in BMO(\mathbb{R}^{n})$ and $m=0,1,\cdots$. Then $m$-th commutators generated by fractional singular integrals $I_{\alpha}$ and $b$ are defined by
$$I^{ m}_{\alpha, b}f(x)=\int_{\mathbb{R}^{n}}\frac{f(y)}{|y-x|^{n-\alpha}}(b(x)-b(y))^{m}dy.$$
and
$$I^{+, m}_{\alpha, b}f(x)=\int_{\mathbb{R}^{n}}\frac{f(y)}{|y-x|^{n-\alpha}}|b(x)-b(y)|^{m}dy.$$
\begin{lemma}\,\,\,\,\cite{ST}\,\,
Let $p,q,\alpha$ and $w$ as in Lemma $2.4$. Then if $b\in BMO(\mathbb{R}^{n})$, there exists $C>0$ such that
both $I^{m}_{\alpha, b}$ and $I^{+,m}_{\alpha, b}$ map $L^{p}(w^{p})$ into $L^{q}(w^{q})$.
\end{lemma}

The high order commutators of degree $m(m\in Z^{+})$ generated by $T_{\alpha}$ and a $BMO(\mathbb{R}^{n})$ function $b$ are defined by
$$(T_{\alpha,b}^{m}f)(x)=\int_{\mathbb{R}^{n}}e^{iP(x,y)}K_{\alpha}(x,y)(b(x)-b(y))^{m}f(y)dy.$$
Lu \cite{L1} first study the $L^{p}$-boundedness for $T_{0,b}^{1}$ when $b\in BMO(\mathbb{R}^{n})$ while the first author of this paper set up the weighted version of the $L^{p}$-boundedness for $T_{b}$ in \cite{Sh}.
In 2005, Wu \cite{W} established the norm estimates for $T_{\alpha,b}^{m}$ when $K_{\alpha}$ is rough kernel. In this section, we will prove the weighted boundedness of $T_{\alpha,b}^{m}$, which is strongly based on Theorem 1.1 and Theorem 1.2.
The method we use here has a root in \cite{DL}.

\begin{theorem} Let $m\in Z^{+},$ $1<p<\infty$, $K_{\alpha}(x,y),P(x,y)$  be the same as in Theorem $1.1$. Then
for $w\in A_{p}$, $b\in BMO(\mathbb{R}^{n})$, there exists constant $C(\alpha)>0, 0<\delta<1$ such that
$$\|T^{m}_{\alpha,b}f\|_{L^{p}(w)}\leq C\|f\|_{L^{p}(w)},$$
where $0<\alpha<(1+\frac{l}{k})\delta$.
\end{theorem}

\begin{theorem} Let $p,m,K_{\alpha},P,b$ be the same as Theorem 3.4 and  $1<p<\frac{n}{\alpha}, \frac{1}{p}-\frac{1}{q}=\frac{\alpha}{n}, 0<\alpha<n$. Then
for $w\in A_{(p,q)}$, there exists constant $C>0$
such that
$$\|T^{m}_{\alpha,b}f\|_{L^{q}(w^{q})}\leq C\|f\|_{L^{p}(w^{p})}.$$
\end{theorem}

{\it Proof of Theorem $3.3$.}
 We carry out the proof of Theorem 3.3 by induction.
First, we consider $m=1$. By Theorem 1.1, $T_{\alpha}$ is bounded on $L^{p}(w)$ with $1<p<\infty$ and $w\in A_{p}$.
Take $\lambda=p$ and $b\in BMO(\mathbb{R}^{n})$. Without loss of generality, we may assume that $\|b\|_{BMO}<\eta$. By Lemma 3.1(1), we have $e^{pb}\in A_{p}$. On the other hand, for every $\theta\in [0,2\pi]$, $b\cos\theta\in BMO(\mathbb{R}^{n})$ and $\|b\cos\theta\|_{BMO(\mathbb{R}^{n})}\leq \|b\|_{BMO(\mathbb{R}^{n})}<\eta$. Thus, $e^{pb\cos\theta}\in A_{p}.$ Now, for $z\in \mathbb{C}$, $g(z)=e^{z(b(x)-b(y))}$ is analytic on $\mathbb{C}$. Thus by the Cauchy integral formula we get
$$b(x)-b(y)=g'(0)=\frac{1}{2\pi i}\int_{|z|=1}\frac{g(z)}{|z|^{2}}dz=\frac{1}{2\pi}\int_{0}^{2\pi}e^{e^{i\theta}[b(x)-b(y)]e^{-i\theta}}d\theta.\eqno(3.1)$$
Since $w\in A_{p}$, by Lemma 2.6(1) we know that there exists $\varepsilon>0$ such that $w^{1+\varepsilon}\in A_{p}$, so
$$\|T_{\alpha}f\|_{L^{p}(w^{1+\varepsilon})}\leq C\|f\|_{L^{p}(w^{1+\varepsilon})}.\eqno(3.2)$$
Let $\lambda=\frac{p(1+\varepsilon)}{\varepsilon}$. By Lemma 3.1(1), we know that there is $\eta>0$ such that $e^{pb(1+\varepsilon)/\varepsilon}\in A_{p}$ whenever $\|b\|_{BMO}<\eta$. Thus for every $\theta\in [0,2\pi]$, we have $e^{pb(1+\varepsilon)\cos\theta/\varepsilon}\in A_{p}$
still holds. By the weighted boundedness of $T_{\alpha}$ we have
$$\|T_{\alpha}f\|_{L^{p}(e^{pb(1+\varepsilon)\cos\theta/\varepsilon})}\leq C\|f\|_{L^{p}(e^{pb(1+\varepsilon)\cos\theta/\varepsilon})}.\eqno(3.3)$$
Applying Lemma 2.7 to (3.2) and (3.3) we have
$$\|T_{\alpha}f\|_{L^{p}(we^{pb\cos\theta})}\leq C\|f\|_{L^{p}(we^{pb\cos\theta})}.\eqno(3.4)$$
Now for $\theta\in [0,2\pi]$, denote $h_{\theta}(x)=f(x)e^{-b(x)e^{i\theta}}$. Then by $f\in L^{p}(w)$, we have
$h_{\theta}\in L^{p}(we^{pb\cos\theta})$ and
$$\|h_{\theta}\|_{L^{p}(we^{pb\cos\theta})}=\|f\|_{L^{p}(w)}.\eqno(3.5)$$
It follows from (3.1) that
\begin{eqnarray*}
T_{\alpha,b}f(x)&=&\int e^{iP(x,y)}K_{\alpha}(x,y)\left(\frac{1}{2\pi}\int_{0}^{2\pi}e^{e^{i\theta}(b(x)-b(y))e^{-i\theta}}d\theta\right)f(y)dy\\
&=&\frac{1}{2\pi}\int_{0}^{2\pi}T_{\alpha}(h_{\theta})(x)e^{e^{i\theta}b(x)}e^{-i\theta}d\theta .
\end{eqnarray*}
Thus by Minkowski's inequality and (3.4) as well as (3.5), we have
$$\|T_{\alpha,b}f\|_{L^{p}(w)}\leq \frac{1}{2\pi}\int_{0}^{2\pi}\|T_{\alpha}(h_{\theta})\|_{L^{p}({we^{pb\cos\theta}})}d\theta\leq C\|f\|_{L^{p}(w)}.$$

 We now assume that Theorem 3.3 holds for $m-1$, i.e. for any $f\in L^{p}(w)$ with $w\in A_{p}$,
$$\left\|T_{\alpha,b}^{m-1}f\right\|_{L^{p}(w)}\leq C\|f\|_{L^{p}(w)}.\eqno(3.6)$$

Using Lemma 2.6(1), there exists $\varepsilon>0$ such that for any $f\in L^{p}(w^{1+\varepsilon})$
$$\left\|T_{\alpha,b}^{k-1}f\right\|_{L^{p}(w^{1+\varepsilon})}\leq C\|f\|_{L^{p}(w^{1+\varepsilon})}.\eqno(3.7)$$
Take $\lambda=\frac{p(1+\varepsilon)}{\varepsilon}$. By Lemma 3.1(1), we know that there is $\eta>0$ such that $e^{pb(1+\varepsilon)/\varepsilon}\in A_{p}$ whenever $\|b\|_{BMO(\mathbb{R}^{n})}<\eta$. On the other hand, for every $\theta\in [0,2\pi]$, $b\cos\theta\in BMO(\mathbb{R}^{n})$ and $\|b\cos\theta\|_{BMO(\mathbb{R}^{n})}\leq \|b\|_{BMO(\mathbb{R}^{n})}<\eta$. Thus, $$e^{pb(1+\varepsilon)\cos\theta / \varepsilon}\in A_{p}.$$ Therefore, by (3.6), we see
that for every $\theta\in [0,2\pi]$ and $f\in L^{p}(e^{pb(1+\varepsilon)\cos\theta / \varepsilon})$, we have
$$\left\|T_{\alpha,b}^{m-1}f\right\|_{L^{p}(e^{pb(1+\varepsilon)\cos\theta / \varepsilon})}\leq C\|f\|_{L^{p}(e^{pb(1+\varepsilon)\cos\theta / \varepsilon})}.\eqno(3.8)$$
Applying Lemma 2.7 to (3.7) and (3.8), we see
that for every $\theta\in [0,2\pi]$ and $f\in L^{p}(we^{pb\cos\theta})$. We have
$$\left\|T_{\alpha,b}^{m-1}f\right\|_{L^{p}(we^{pb\cos\theta})}\leq C\|f\|_{L^{p}(we^{pb\cos\theta})}.\eqno(3.9)$$
Moreover, set $h_{\theta}(x)=f(x)e^{-b(x)e^{i\theta}}$. Then by $f\in L^{p}(w)$, we have for $\theta\in [0,2\pi]$,
$h_{\theta}\in L^{p}(we^{pb\cos\theta})$ and
$$\|h_{\theta}\|_{L^{p}(we^{pb\cos\theta})}=\|f\|_{L^{p}(w)}.\eqno(3.10)$$
For simplicity, we denote $$K_{\alpha,m-1}(x,y)=e^{iP(x,y)}K_{\alpha}(x,y)(b(x)-b(y))^{m-1}\chi_{\{2^{j-1}<|y-x|<2^{j}\}}(x,y).$$
From (3.1), we have
\begin{eqnarray*}
T_{\alpha,b}^{m}f(x)&=&\int K_{\alpha,m-1}(x,y)(b(x)-b(y))f(y)dy\\
&=&\frac{1}{2\pi}\int_{0}^{2\pi}T^{m-1}_{\alpha, b}(h_{\theta})(x)e^{e^{i\theta}b(x)}e^{-i\theta}d\theta.
\end{eqnarray*}
By (3.9) and (3.10), we get
$$\left\|T_{\alpha,b}^{m}f\right\|_{L^{p}(w)}\leq C\|f\|_{L^{p}(w)}.$$
We complete the proof of Theorem 3.3.
 \qed

{\it Proof of Theorem $3.4$.}
As a result of the similar propositions of $A_{p}$ and $A_{(p,q)}$ in Lemma 2.6 as well as Lemma 3.1(2), the proof of Theorem 3.4 can be handled in much the same way as that of Theorem 3.3, which we only need to do a slight modifications by replacing $A_{p}$ by $A_{(p,q)}$.  However, by the idea similar to that of in the proof of Theorem 1.2, we can also easily obtain Theorem 3.4 by Lemma 3.2
and the following observation
$$
\left|T_{\alpha,b}f(x)\right|\leq \int\frac{|f(y)(b(x)-b(y))|}{|x-y|^{n-\alpha}}=I_{\alpha,b}(|f|)(x).
$$
 \qed

\end{document}